\documentclass[12pt]{article}
\usepackage{amssymb, amsmath, amsthm}
\usepackage[all]{xy}

\newcommand{\mr}[1]{\mathrm{#1}}
\newcommand{\mf}[1]{\mathfrak{#1}}
\newcommand{\mc}[1]{\mathcal{#1}}

\newcommand{\Z}{{\bf Z}}
\newcommand{\Q}{{\bf Q}}
\newcommand{\zp}{{\bf Z}_p}

 \DeclareMathOperator{\Gal}{Gal}
 
\DeclareMathOperator{\coker}{coker}

\newtheorem{theorem}{Theorem}[section]
\newtheorem{proposition}[theorem]{Proposition}
\newtheorem{lemma}[theorem]{Lemma}
\newtheorem{corollary}[theorem]{Corollary}

\newtheorem{conjecture}[theorem]{Conjecture}

\theoremstyle{definition}

\theoremstyle{remark}

\renewcommand{\baselinestretch}{1.2}

\begin{document}

\title{On Galois groups of unramified pro-$p$ extensions}
\author{Romyar T. Sharifi}
\date{}
\maketitle

\begin{abstract}
	Let $p$ be an odd prime satisfying Vandiver's conjecture.
	We consider two objects, the Galois group $X$ 
	of the maximal unramified abelian pro-$p$
	extension of the compositum of all $\zp$-extensions of $\Q(\mu_p)$ and
	the Galois group $\mf{G}$ of the maximal unramified pro-$p$ extension of 
	$\Q(\mu_{p^{\infty}})$.  We give a lower bound for the height of the
	annihilator of $X$ as an Iwasawa module.  Under some mild assumptions
	on Bernoulli numbers, we provide a necessary and sufficient
	condition for $\mf{G}$ to be abelian.  The bound and the condition in 
	the two results are given in terms of special values of a cup product
	pairing on cyclotomic $p$-units.  We obtain in particular that, for $p <
	1000$, Greenberg's conjecture that $X$ is pseudo-null holds
	and $\mf{G}$ is in fact abelian.
\end{abstract}

\section{Introduction}

Let $L$ be a possibly infinite algebraic extension of $\Q$, and fix a prime $p$.  We shall say that an algebraic extension of a number field is unramified if all of its finite subextensions are unramified at all places.  The following are two frequently recurring questions in algebraic number theory.

\begin{enumerate}
\item[1.] 
What is the structure of the Galois group $X_L$ of the maximal unramified abelian pro-$p$ extension of $L$?  
\item[2.] 
What is the structure of the Galois group $\mf{G}_L$ of the maximal unramified pro-$p$ extension of $L$?
\end{enumerate}

If $L$ is a number field, one sort of answer to the first question is found in class field theory.  The maximal unramified abelian pro-$p$ extension of $L$ is finite and is known as the $p$-Hilbert class field, and $X_L$ is isomorphic to the $p$-part $A_L$ of the class group of $L$.  As for the second question, the maximal unramified pro-$p$ extension of $L$ can be infinite, as Golod and Shafarevich demonstrated the existence of infinite Hilbert $p$-class field towers.

Consider next the cyclotomic $\zp$-extension $K$ of a number field $F$.  We let $\Gamma = \Gal(K/F)$ and set $\Lambda = \zp[[\Gamma]]$.  Iwasawa showed that the group $X_K$ is always finitely generated and torsion over $\Lambda$.  
Beyond that, if $F$ is an abelian extension of $\Q$, the now-proven ``main conjecture of Iwasawa theory" states the characteristic ideal of an odd eigenspace of $X_K$ is determined by the $p$-adic $L$-function of a related even character.  Wiles proved a similar result for abelian characters over
totally real number fields.

In most instances, little is known about the structure of $\mf{G}_K$ that cannot
be obtained from the structure of $X_K$, its maximal abelian quotient.  We consider the question of whether or not $\mf{G}_K$ is abelian in the fundamental case that $F = \Q(\mu_p)$ for an odd prime $p$.  

\begin{theorem}
	For $p < 1000$, the group $\mf{G}_K$ is abelian.  For $p = 1217$, 
	$7069$, and $9829$, it is not.
\end{theorem}

For the other 241 primes $p <$ 25,000 such that the $\zp$-rank of $X_K$ is at least $2$, a conjecture of McCallum and the author's \cite[Conjecture 5.3]{mcs} implies that $\mf{G}_K$ is a free abelian pro-$p$ group, so $\mf{G}_K \cong X_K$.  For the general result, see Theorem \ref{GKresult}.  Until recently, it was thought to be proven that $\mf{G}_K$ is free pro-$p$ under Vandiver's conjecture that $p$ does not divide the class number of $\Q(\mu_p)^+$ \cite{nqd,wingberg}, e.g., for $p <$ 12,000,000.  Now, it seems quite likely that $\mf{G}_K$ is never free pro-$p$ unless it is trivial or isomorphic to $\zp$.

For larger fields $L$, one can ask about the ``size" of $X_L$.  Early on, Greenberg proposed a conjecture in the abelian setting (see \cite[Conjecture 3.5]{greenberg}).  Let $\tilde{F}$ denote the compositum of all $\zp$-extensions of $F$, let $\tilde{\Gamma} = \Gal(\tilde{F}/F)$, and let $\tilde{\Lambda} = \zp[[\tilde{\Gamma}]]$.  It follows from the analogous statement for $K$ that $X_{\tilde{F}}$ is finitely generated and torsion as a $\tilde{\Lambda}$-module.  To get slightly more information, we can pass to a family of $\tilde{\Lambda}$-modules that is in a definite sense one step smaller than the family of $\tilde{\Lambda}$-torsion modules.  We say that a finitely generated $\tilde{\Lambda}$-module is pseudo-null if its annihilator has height at least $2$.

\begin{conjecture}[Greenberg] \label{greenconj}
	The Galois group $X_{\tilde{F}}$ of the maximal unramified
	abelian pro-$p$ extension of $\tilde{F}$ is pseudo-null as a 
	$\tilde{\Lambda}$-module.
\end{conjecture}

In Theorem \ref{bettergreen}, we give a lower bound on the height of the annihilator of $X_{\tilde{F}}$ for $p$ satisfying Vandiver's conjecture.  This implies the following weaker statement.

\begin{theorem} \label{green1000}
	For $p < 1000$, Greenberg's conjecture holds for $p$ and $F = \Q(\mu_p)$.
\end{theorem}

We can ask a similar question for the class of strongly admissible $p$-adic
Lie extensions $L$ of $F$.  As in \cite{hs}, we say that $L/F$ is strongly admissible if it is ramified at only finitely many places of $F$, $L$ contains $K$, and ${Q} = \Gal(L/F)$ is pro-$p$, has dimension at least $2$, and contains no elements of order $p$.  

We define a finitely generated $\Lambda({Q}) = \zp[[{Q}]]$-module 
$M$ 
to be pseudo-null if 
$$
	\mr{Ext}^i_{\Lambda({Q})}(M			
	,\Lambda({Q})) = 0
$$ 
for $i = 0,1$.  Note that ${Q}$ contains a normal subgroup $G = \Gal(L/K)$ with quotient $\Gamma$.  
If $M$ 
is finitely generated over $\Lambda(G)$, 
then 
$M$ 
is $\Lambda({Q})$-pseudo-null if and only if 
$M$
is $\Lambda(G)$-torsion \cite[Lemma 3.1]{hs}. 
Under Iwasawa's conjecture on the triviality of the $\mu$-invariant of $X_K$, the finite generation of $X_L$ over $\Lambda(G)$ would always hold for $L/F$ strongly admissible \cite[Lemma 3.4]{hs}.  
The question is then to determine the $L$ as above for which $X_L$ is $\Lambda(G)$-torsion.  In particular, we use this in showing the following.

\begin{theorem} \label{pnexample}
	Let $F = \Q(\mu_p)$ for some $p < 1000$, and suppose that $L$ is a strongly
	admissible $p$-adic Lie extension of $F$ 
	that contains a $\zp$-extension 
	of $K$ that is unramified outside $p$ and contains a $p$th root of $p$.
	Then $X_L$ is $\Lambda({Q})$-pseudo-null.
\end{theorem}

\section{Growth of Iwasawa modules}

For simplicity of the description, let us assume that either $p$ is odd or our
number field $F$ is purely imaginary.  We use $S$ to denote a set of primes of $F$ containing those above $p$, and for any algebraic extension $E$ of $F$, we let $S_E$ denote the set of primes of $E$ above those in $S$.
We use $\mf{X}_{E,S}$ to denote the Galois group of the maximal abelian pro-$p$ unramified outside $S_E$ extension of $E$.  For $v \in S_E$, we let $G_{E_v}$ (resp., $I_{E_v}$) denote the absolute Galois group (resp., inertia subgroup) of the completion $E_v$, and we
let $\mc{V}_E$ (resp., $\mc{U}_E$) denote the inverse limit under norm maps of the $p$-completions of the unit
groups (resp., $p$-unit groups) of number fields in $E$.  
Let 
$$
	\mc{W}_E = 
	\lim_{\substack{\leftarrow\\F' \subset E}} 
	\bigoplus_{v \in S_{F'}} H_1(I_{F'_v},\zp)_{G_{F'_v}},
$$
where $F'$ runs over the finite extensions of $F$ in $E$.
We have the following proposition.

\begin{proposition} \label{torsion}
	Let $E$ be an algebraic extension of $F$, and let $L
	\subsetneq M$ be pro-$p$ $p$-adic Lie extensions of $E$ such that
	$G = \Gal(M/E)$ has no $p$-torsion.  Let $H = \Gal(M/L)$.  
	Let $S$ be the set of primes of $F$ consisting of
	those above $p$ and those that ramify in $M/F$.  
	\begin{enumerate}
		\item[a.] Suppose that the set of primes of $S_L$ 
		that ramify in $M/L$ is finite and 
		that there exists a $\zp$-extension $E'$ of 
		$E$ contained in $L$ such that $X_L$ is finitely generated over 
		$\Lambda(\Gal(L/E'))$. Then $X_M$ is finitely generated and
		torsion over $\Lambda(G)$.
		\item[b.]   Suppose that there exists exactly one prime in $S_L$ 
		that ramifies in $M/L$ and 
		that $M/L$ is totally ramified at that prime. Then 
		$(X_M)_H \cong X_L$ via restriction.
	\end{enumerate}
\end{proposition}

\begin{proof}
	We consider the natural commutative diagram
	$$
		\SelectTips{eu}{} \xymatrix@R=20pt@C=15pt{
		\mc{V}_M \ar[r] \ar[d]^a & \mc{W}_M \ar[r] \ar[d]^b &
		\mf{X}_{M,S} \ar[r] \ar[d]^c & X_M \ar[r] \ar[d]^d & 0 \\
		\mc{V}_L \ar[r] & \mc{W}_L \ar[r] & \mf{X}_{L,S} \ar[r] & X_L \ar[r] & 0.
		}
	$$
	Let $b_H$, $c_H$ and $d_H$ denote the maps 
	induced by $b$, $c$ and $d$, respectively, on $H$-coinvariants of 
	the domains.
	A straightforward diagram chase yields an isomorphism
	\begin{equation} \label{dcresult}
		\coker(\ker c_H \to \ker d_H) \cong
		\frac{\ker(\coker b \to \coker c)}{\mr{image}(\coker a \to \coker b)}.
	\end{equation}
	By the Hochschild-Serre spectral sequence, we have that $\ker c_H$
	is isomorphic to a quotient of $H_2(H,\zp)$, which is finitely generated over
	$\zp$.
	Since we have assumed that the set of primes of 		
	$S_L$ at which $M/L$ ramifies is finite, we have that
	$$
		\coker b \cong \bigoplus_{v \in S_L} J_v,
	$$
	where $J_v$ is the inertia subgroup of 
	the abelianization of the decomposition group at $v$ in $H$.
	Thus $\coker b$ is finitely generated over $\zp$.
	It follows that $\ker d_H$ is also finitely generated over $\zp$.
	Since, by assumption, $X_L$ is finitely generated over $\Lambda(\Gal(L/E'))$,
	the same must hold for $(X_M)_H$.  By \cite[Section 3]{bh}, we then have that $X_M$ is finitely 
	generated over $\Lambda(\Gal(M/E'))$, and by \cite[Proposition 2.3]{cfksv}, we 
	may conclude that $X_M$ is a torsion $\Lambda(G)$-module,
	proving part a.
	
	Now, suppose that the assumption of part b holds.  
	Let $w$ denote the unique prime of $L$ at which $M/L$
	is totally ramified, and let $I_w$ denote the inertia group of 
	$w$ in $H$.  The existence of such a prime forces
	$\coker d = 0$ automatically.
	Since $I_w = H$, we also have that the map $\coker b \to \coker c$ is 
	an isomorphism.  Furthermore, 
	$\ker b_H$ is canonically isomorphic to a quotient of $H_2(I_w,\zp)$,
	and $H_2(I_w,\zp) \to H_2(H,\zp)$ is an isomorphism since $I_w = H$.
	Thus the map $\ker b_H \to \ker c_H$ is surjective, and this forces the map
	$\ker c_H \to \ker d_H$ to be zero.  By \eqref{dcresult}, we therefore
	have that $\ker d_H = 0$, proving part b.
\end{proof}

Next, consider the following lemma from commutative algebra.

\begin{lemma} \label{heights}
	Let $n \ge m \ge 0$, and consider the power series ring 
	$R_n = \zp[[T_1,\ldots,T_n]]$ in $n$ independent variables.  
	Suppose that $A$ is an $R_n$-module
	that is finitely generated over $R_m$.  Let $h$ denote the height of the
	annihilator of $A$ in $R_m$.  Then the height of the
	annihilator of $A$ in $R_n$ is at least $h+n-m$.
\end{lemma}

\begin{proof}
	Let $I_m$ (resp., $I_n$) denote the annihilator of $A$
	in $R_m$ (resp., $R_n$).  Suppose that $A$ is nonzero and 
	generated over $R_m$ by a given finite set of elements.  
	For $d$ with $m < d \le n$, 
	we may consider any characteristic polynomial $f_d \in
	R_m[X]$ of $T_d$ acting on $A$ with respect to this generating set.
	Now, if $\mf{p}$ is a prime ideal containing $I_n$, then it contains $I_m$.
	Let $\mf{q}_1 \subsetneq \ldots \subsetneq \mf{q}_h$ 
	be a chain of distinct nonzero primes of 
	$R_m$ contained in $\mf{p} \cap R_m$.
	For $d$ as above, set $i = d+h-m$, and 
	define $\mf{q}_i = \mf{p} \cap R_d$.  Then $f_d(T_d) \in \mf{q}_i$,
	but $f_d(T_d) \notin \mf{q}_{i-1}R_d$ since $f_d$ is monic.
	Let $\mf{p}_i = \mf{q}_iR_n$ for all $i$ with $1 \le i \le h+n-m$.  
	Then the $\mf{p}_i$ are prime, distinct, and 
	contained in $\mf{p}$, so the height of $I_n$ is
	at least $h+n-m$.
\end{proof}

The following corollary will be useful to us later.

\begin{corollary} \label{growthcor}
	Let $L \subseteq M$ be abelian extensions of $F$ containing 
	the cyclotomic $\zp$-extension	 $K$ of $F$ with 	
	$\Gal(M/F) \cong \zp^n$ and $\Gal(L/F) \cong \zp^m$ for some $n
	\ge m \ge 1$.  If $X_L$ is finitely generated over $\zp$,
	then the height of the annihilator of $X_M$ is at least $m$
	as a $\Lambda(\Gal(M/F))$-module.
\end{corollary}

\begin{proof}
	Since $M/F$ is necessarily unramified outside $p$ and there are only
	finitely many primes above $p$ in $L$, the assumptions of Proposition
	\ref{torsion}a are satisfied, so $X_M$ is a finitely generated,
	torsion $\Lambda(\Gamma \times \Gal(M/L))$-module.  
	The result is then immediate from Lemma \ref{heights}.
\end{proof}

\section{Iwasawa modules over Kummer extensions}

For the rest of this paper, we will focus on the case that $F = \Q(\mu_p)$
for an odd prime $p$.  The key to the results in this paper is the use of 
cup products to control the growth of the Galois groups
of the maximal unramified abelian pro-$p$ extensions of
Kummer extensions of $K = \Q(\mu_{p^{\infty}})$.
These cup products were first studied in \cite{mcs} and later in \cite{me-paireis}.  
Let us briefly review the material we will need.

Let $\mc{B}_F$ denote the $p$-completion of the
group of elements of $F^{\times}$ whose $p$th roots generate unramified outside $p$ extensions of $F$.  Given an element of the $p$-completion $\mc{K}$ of $K^{\times}$, we may also speak of the extension of $K$ its $p$-power roots generate, and we let $\mc{B}_K$ denote the subgroup of those elements whose $p$-power roots generate unramified outside $p$ extensions.  (Although $\mc{B}_F$ is not a subgroup of $\mc{B}_K$, there is a canonical map from the $p$-completion of $F^{\times}$ to $\mc{K}$ with kernel $\mu_p$ through which we will consider elements of $\mc{B}_F$ as elements of
$\mc{K}$ by abuse of notation.)

Let $S$ now denote the set consisting of the unique prime above $p$ in $F$, and for any algebraic extension $E$ of $F$, let $G_{E,S}$ denote the Galois group of the maximal unramified outside $S_E$ extension of $E$.  The cup product on $H^1(G_{F,S},\mu_p)$ induces a pairing
$$
	(\,\cdot\,,\,\cdot\,)_{p,F,S} \colon \mc{B}_F \times \mc{B}_F \to A_F \otimes \mu_p,
$$
where $A_F$ again denotes the $p$-part of the class group of $F$,
and the inverse limit of similar pairings arising from the cup products on the groups $H^1(G_{\Q(\mu_{p^n}),S},\mu_{p^n})$ induces a pairing
$$
	(\,\cdot\,,\,\cdot\,)_{K,S} 
	\colon \mc{B}_K \times \mc{U}_K \to X_K(1)
$$
(see, for instance, \cite[Section 4]{me-paireis}).

Fix a primitive $p$th root of unity $\zeta_p$ in $F$. 
Let $\mc{E}_F$ denote the $p$-completion of the $p$-units in $F$.
Let $\mc{C}_F$ denote the subgroup of $\mc{E}_F$ consisting
of the $p$-completion of the cyclotomic $p$-units, i.e., those generated by
the $1-\zeta_p^i$ with $1 \le i \le p-1$.  
For a finite extension $E$ of $F$, let $A_{E/F}$ denote the image of
the norm map $A_E \to A_F$.  For an algebraic extension $L/K$, let
$X_{L/K}$ denote the image of restriction $X_L \to X_K$,
and let $Y_L$ denote the maximal quotient of $X_L$ in which
all primes above $p$ split completely.  
For a $p$-adic Lie group $G$, let $I_G$ denote the augmentation ideal
in $\Lambda(G)$.  It will be useful to recall the following result.

\begin{lemma} \label{augtriv}
	Let $b \in \mc{B}_K$, and suppose that its $p$-power roots
	generate a $\zp$-extension
	$L$ of $K$ with Galois group $G$.  We have that
	\begin{equation} \label{surj1}	
		(b,\mc{U}_K)_{K,S} = X_{L/K}(1)
	\end{equation}
	if and only if $I_G Y_L = 0$.  If $ba^{-1} \in 
	\mc{K}^p$ for some $a \in \mc{B}_F$, then \eqref{surj1} follows from
	\begin{equation} \label{surjectivity}
		(a,\mc{C}_F)_{p,F,S} = A_{F(a^{1/p})/F} \otimes \mu_p.
	\end{equation}
	Moreover, if Vandiver's conjecture holds at $p$, then \eqref{surj1} implies
	\eqref{surjectivity} for some $a$ with $ba^{-1} \in \mc{K}^p$.
\end{lemma}

\begin{proof}
	We remark that the containment of the pairing values in $X_{L/K}(1)$ (resp.,
	in $A_{F(a^{1/p})/F} \otimes \mu_p$) in the statement
	holds for any $b \in \mc{B}_K$ (resp., $a \in \mc{B}_F$), 
	as follows for example from \cite[Theorem 2.4]{mcs}.
	The first part 
	is a direct consequence of \cite[Theorem 4.3]{me-paireis},
	noting that there is a unique prime over $p$ in $K$.  
	Now, we know that the image of $(b,\mc{U}_K)_{K,S}$ under the natural
	map $X_K(1) \to A_F \otimes \mu_p$
	contains $(a,\mc{C}_F)_{p,F,S}$ and is equal to it under Vandiver's conjecture at $p$,
	since every element of $\mc{C}_F$ is a universal norm from $K$
	and these are all the universal norms under Vandiver.  Since we have as well that
	the image of $X_{L/K}(1)$ is contained in	
	$A_{F(a^{1/p})/F} \otimes \mu_p$ and is equal to it when $a$ is chosen properly 
	(multiplying $a$ by an element of $\mu_p$ to make 
	$F(a^{1/p})/F$ unramified if possible),
	the second statement follows.
\end{proof}

We next recall \cite[Corollary 5.9]{me-paireis}.

\begin{theorem} \label{pairp}
	For $p < 1000$, we have $(p,\mc{C}_F)_{p,F,S} = A_F \otimes \mu_p$.
\end{theorem}

We now obtain the following useful result.

\begin{proposition} \label{surjtors}
	Suppose that $a \in \mc{B}_F$ is
	such that \eqref{surjectivity} holds.
	Let $L$ be a $\zp$-extension of $K$ with Galois group $G$ 
	that is unramified outside $p$
	but not unramified and contains
	a $p$th root of $a$.  Then $X_L$ is
	finitely generated over $\zp$, 
	and we have $Y_L \cong X_{L/K}$ via restriction.
\end{proposition}

\begin{proof}
	The field $L$ is defined over $K$
	by the $p$-power roots of an element $b \in \mc{B}_K$  
	such that $ba^{-1} \in \mc{B}_K^p$. 
	We have $I_G Y_L = 0$ by Lemma \ref{augtriv}.
	It follows that $Y_L \cong (Y_L)_G$, and $(Y_L)_G \cong X_{L/K}$ is
	easily seen since $L/K$ is a $\zp$-extension with a unique prime over
	$p$ in $K$.  In addition, we know that $(Y_L)_G = (X_L)_G$, so $I_G X_L$
	is generated by the decomposition groups above $p$ in $X_L$.  	
	Such decomposition groups are attached to the primes above $p$ in $L$, of
	which there are finitely many since $L/K$ is not unramified, and any
	such decomposition group is a quotient of $\zp$.  Therefore, $X_L$
	is a finitely generated $\zp$-module.
\end{proof}
 
Noting the remarks preceding it in the introduction,
Theorem \ref{pnexample} now follows from Theorem \ref{pairp}, Proposition \ref{surjtors}, and 
Proposition \ref{torsion}a.

\section{Iwasawa modules over multiple $\zp$-extensions}

In this section, we turn to the study of the structure of the Galois group $X_{\tilde{F}}$ of the maximal unramified abelian pro-$p$ extension of the compositum $\tilde{F}$ of all $\zp$-extensions of the $p$th cyclotomic field $F$.  Proposition \ref{surjtors} allows us to immediately give a criterion for the verification of Greenberg's pseudo-nullity conjecture (Conjecture \ref{greenconj}) for $F$.

\begin{proposition}
	Suppose that there exists $a \in \mc{C}_F$ 
	satisfying \eqref{surjectivity}.
	Then Greenberg's conjecture holds for $F$ and $p$.  
\end{proposition}

\begin{proof}
	Since $(\zeta,\mc{C}_F)_{p,F,S} = 0$
	for $\zeta \in \mu_p$, we may assume without loss of generality 
	that $a$ is not a root of unity times a $p$th power in $\mc{E}_F$.  
	Since $a \in \mc{C}_F$, there exists a $\zp$-extension $E$ of $F$
	containing a $p$th root of $a$.  By assumption, $E$ is not contained in $K$.
	Setting $L = EK$, we have that
	$X_L$ is finitely generated over $\zp$ by Proposition \ref{surjtors}, and
	since $L \subseteq \tilde{F}$, that Greenberg's conjecture holds
	by Proposition \ref{torsion}a, noting \cite[Lemma 3.4]{hs} (or 
	\cite[Proposition 5.4]{venjakob}).
\end{proof}

In particular, as stated in Theorem \ref{green1000}, Greenberg's conjecture holds for $\Q(\mu_p)$ for $p < 1000$.   We remark that Theorem \ref{green1000} was already known for those $p < 1000$ with $A_F$ cyclic by \cite[Corollary 10.5]{mcs} and Theorem \ref{pairp}. 
An inductive version of the above argument is due to Greenberg in the case that $X_K \cong \zp$, but unlike the argument of \cite{mcs}, it works easily without this restriction.  We thank Ralph Greenberg for communicating his argument to us.

One may ask if the pseudo-nullity of $X_{\tilde{F}}$ is the best one can do.  That is, can one give a stronger lower bound on the height of the annihilator of $X_{\tilde{F}}$?  In fact, the answer is yes, as we now show.  

Let us say that $k$ is irregular for $p$ if $(p,k)$ is an irregular pair, i.e., $k$ is positive and even, $k \le p-3$, and $p$ divides the $k$th Bernoulli number $B_k$.  
Let $\Delta = \Gal(F/\Q)$.  For a $\zp[\Delta]$-module $A$
and $i \in \Z$, let $A^{(i)}$ denote the $\omega^i$-eigenspace of $A$,
where $\omega$ denotes the Teichm\"uller character.  Then $(p,k)$ is irregular 
if and only if $A_F^{(1-k)} \neq 0$.  
For any odd integer $i$, let $\eta_i$ denote the projection of $(1-\zeta_p)^{p-1}$ to $\mc{C}_F^{(1-i)}$.

\begin{theorem} \label{bettergreen}
	Suppose that Vandiver's conjecture holds at $p$.
	Consider the following subsets of $\Z/(p-1)\Z$:
	$$
		R = \{ k \mid (p,k) \mr{\ irregular} \}
	$$  			
	and 
	$$
		I = \{ i \mid i \mr{\ odd}, (\eta_i,\eta_{k-i})_{p,F,S}
		\neq 0 \mr{\ for\ all\ } k \mr{\ irregular\ for\ } p \}.
	$$   
	The height of the annihilator of $X_{\tilde{F}}$ as a 
	$\tilde{\Lambda}$-module is at least one more than
	the maximal number of disjoint translates $i + R$ with $i \in I$.
\end{theorem}

\begin{proof}
	Since we have assumed Vandiver's conjecture, $A_F^{(1-k)}$ is
	cyclic for all irregular $k$ for $p$, and all other eigenspaces of 
	$A_F$ are trivial.  Therefore, the condition that $i \in I$ 
	is equivalent to that of \eqref{surjectivity} holding for $a = \eta_i$,
	along with the extension of $F$ defined by a $p$th root of $\eta_i$ being
	totally ramified above $p$ (i.e., 
	$i \neq p-k'$ for all $k'$ irregular for $p$ \cite[Section 5]{mcs}). 
	Let $L_i$ denote the unique $\zp$-extension of $K$ Galois over $\Q$ 
	and abelian over $F$ that contains a $p$th root of $\eta_i$.  
	Then Proposition \ref{surjtors} yields that
	$Y_{L_i} \cong X_K$ via restriction.
	
	Let $i_1, i_2, \ldots, i_d \in I$ be such that the translates
	$i_s + R$ are all disjoint as $s$ runs over $1 \le s \le d$, 
	let $M_s = L_{i_1} \cdots L_{i_s}$ for any such $s$, and set $M_0 = K$.  
	Suppose by induction on $d$ that $Y_{M_{d-1}} 
	\cong X_K$ via restriction.  Again, the assumption that $i \in I$
	implies that $L_i/K$ is totally ramified at $p$.  Since the
	$\Gal(L_i/K)$ have $\Delta$-actions given by distinct powers of
	the Teichm\"uller character, $M_d$ has a unique prime
	above $p$, and that prime is totally ramified over $F$.  
	
	Now set $G = \Gal(M_d/K)$, $H = \Gal(M_d/M_{d-1})$, and $T
	= \Gal(M_d/L_d)$.  By
	Proposition \ref{torsion}b, we know that 
	$(Y_{M_d})_T \cong Y_{L_d}$.
	We therefore have
	\begin{equation} \label{zeroquot}
		0 = (I_H Y_{L_d})_H \cong (I_H (Y_{M_d})_T)_H \cong
		I_H Y_{M_d}/(I_T Y_{M_d} \cap I_H Y_{M_d} + I_H^2 Y_{M_d}).
	\end{equation}
	By assumption, we have
	$$
		I_T Y_{M_d}/(I_T Y_{M_d} \cap I_H Y_{M_d})
		\cong I_T((Y_{M_d})_H) \cong
		I_TY_{M_{d-1}} \cong I_T X_K = 0,
	$$
	so we have
	\begin{equation} \label{eqint}
		I_T Y_{M_d} = I_T Y_{M_d} \cap I_H Y_{M_d}
	\end{equation}
	and therefore
	$$
		I_T Y_{M_d} \subseteq I_H Y_{M_d}.
	$$
	
	Consider for $N = H$ and $N = T$ 
	the natural surjective $\zp[\Delta]$-homomorphisms
	$$
		\pi_N \colon X_K \otimes_{\zp} N \to 
		(I_N Y_{M_d})_G,
	$$
	with
	$$
		\pi_N(x \otimes \sigma) = (\sigma-1)\tilde{x} \pmod{I_GI_N Y_{M_d}},
	$$
	where $\tilde{x} \in Y_{M_d}$ restricts to $x$.
	Since the $\zp[\Delta]$-eigenspaces of
	$X_K \otimes_{\zp} N$ are nontrivial outside of those of the characters
	$\omega^{1-k+i_t}$ with $k$ irregular for $p$ and $t \le d-1$ if $N = T$
	and $t = d$ if $N = H$, we have that $(I_N Y_{M_d})_G$
	is also nontrivial at most in these eigenspaces.  Since the $i_t + R$
	are all disjoint, the canonical map
	$$
		(I_T Y_{M_d})_G \to (I_H Y_{M_d})_G
	$$ 
	is zero, and hence
	$$
		I_T Y_{M_d} \subseteq I_GI_H Y_{M_d} = (I_TI_H + I_H^2)Y_{M_d}.
	$$
	But \eqref{eqint} then forces
	$$
		(I_T Y_{M_d} \cap I_H Y_{M_d}) + I_H^2 Y_{M_d}
		= I_G I_H Y_{M_d}.
	$$
	Given this, \eqref{zeroquot} implies that $(I_H Y_{M_d})_G = 0$,
	that is,
	$$
		Y_{M_d} \cong (Y_{M_d})_H \cong Y_{M_{d-1}},
	$$ 
	and so $Y_{M_d} \cong X_K$ via restriction.  
	
	Since there exists a unique prime over $p$ in $M_d$, the kernel of 
	$X_{M_d} \to Y_{M_d}$ is a quotient of $\zp$, so $X_{M_d}$ is finitely
	generated over $\zp$.
	By Corollary \ref{growthcor}, the annihilator of $X_{\tilde{F}}$ then has
	height at least $d+1$ as a $\tilde{\Lambda}$-module.  
\end{proof}

\begin{corollary}
	Suppose that Vandiver's conjecture holds at $p$, let $r$ be the
	$p$-rank of $A_F$, and let $s$ be the number of
	odd integers $i$ with $1 \le i \le p-2$ such that
	$(\eta_i,\eta_{k-i})_{p,F,S} \neq 0$ for all $k$ irregular for $p$.
	Then the height $j$ of the annihilator of 
	$X_{\tilde{F}}$ as a $\tilde{\Lambda}$-module satisfies 
	$$
		j \ge \frac{s}{r^2-r+1} + 1.
	$$
\end{corollary}

\begin{proof}
	Let $R$ and $I$ be as in Theorem \ref{bettergreen}.
	Suppose that we have $t$ disjoint translates
	of $R$ by elements of $I$. 
	Each translate of $R$ contains $r$ elements and 
	intersects at most $r^2-r+1$ of the $s$ translates of $R$ by elements of 
	$I$.  Therefore, so long at $t(r^2-r+1) < s$, there exists at least
	one other translate by an element of $I$ that does not intersect any of 
	the given translates.  The result then follows from Theorem 
	\ref{bettergreen}.
\end{proof}

We remark that, for $p < 1000$, one has $r \le 3$ and 
$\frac{p-1}{2}-s \in [2,6]$, $[6,8]$, and $[9,12]$ when $r = 1$, $2$, 
and $3$, respectively.

\section{The maximal unramified pro-$p$ extension}

We now turn to the study of the structure of the Galois group $\mf{G}_K$ of the maximal unramified pro-$p$ extension of the field $K$ of all $p$-power roots of unity.

\begin{lemma} \label{keyseq}
	Let $M/L$ be an unramified
	abelian pro-$p$ extension with torsion-free Galois group $H$, and   
	assume that $M$ and $L$ are Galois extensions of $F$. 
	Then there is a canonical exact sequence of $\Lambda(\Gal(L/F))$-modules,
	$$
		H \wedge_{\zp} H \to (X_M)_H \to X_L \to H \to 0.
	$$
\end{lemma}

\begin{proof}
	This arises directly from the obvious exact sequence
	$$
		1 \to \frac{[\mf{G}_L,\mf{G}_L]}{[\mf{G}_M,\mf{G}_L]} \to 
		\frac{\mf{G}_M}{[\mf{G}_M,\mf{G}_L]}
		\to \frac{\mf{G}_L}{[\mf{G}_L,\mf{G}_L]} \to 
		\frac{\mf{G}_L}{\mf{G}_M[\mf{G}_L,\mf{G}_L]} \to 1.
	$$
\end{proof}

Moreover, we note that the sequences in Lemma \ref{keyseq} are natural in
$L$ and $M$ (and thereby $H$) in the obvious sense.
We have the following result on the structure of $\mf{G}_K$.   

\begin{theorem} \label{GKresult}
	Assume that Vandiver's conjecture holds at $p$ and 
	that $X_K^{(1-k)} \cong \zp$ 
	for each irregular $k$ for $p$.
	Assume also that for any irregular $k' > k$ for $p$, 
	we have both that $k + k' \not\equiv 2 \bmod p-1$ and that $j = k$
	whenever $j' > j$ are irregular for $p$ with 
	$j + j' \equiv k + k' \bmod p-1$.
	Then the group $\mf{G}_K$ is abelian if and only if 
	$(\eta_{p-k},\eta_{k+k'-1})_{p,F,S} \neq 0$ for all irregular $k' > k$
	for $p$.
\end{theorem}

\begin{proof}
	Under Vandiver's conjecture, $X_K$ is $p$-torsion free, so $\mf{G}_K$
	is free pro-$p$ and abelian if the $\zp$-rank of $X_K$ is at most $1$.
	Hence, we may assume that the $\zp$-rank of $X_K$ is at least $2$.  
	Then, by assumption, 
	$X_K$ has at least two nontrivial $\zp[\Delta]$-eigenspaces, 
	where $\Delta = \Gal(F/\Q)$. Fix $k$ 
	irregular for $p$, and consider the unramified $\zp$-extension $L/K$ defined 
	by $H = X_K^{(1-k)}$.  
		
	We begin by showing that $X_L \cong (X_L)_H$.  This is equivalent to
	showing that $(I_H X_L)_H$ is trivial.  Note first that  $(X_L)_H \cong (Y_L)_H$.
	Furthermore, we claim that
	\begin{equation} \label{samequot}
		(I_H X_L)_H \cong (I_H Y_L)_H.
	\end{equation}
	The image of a decomposition group above $p$ in $X_L$ 
	must lie in the trivial
	$\zp[\Delta]$-eigenspace of $(I_H X_L)_H$.  On the other hand, by 
	the procyclicity of the eigenspaces of $X_K$ and
	Vandiver's conjecture,
	$(I_H X_L)_H$ can only be nontrivial in its 
	$(\omega^{2-k-k'})$-eigenspaces for $k'$ irregular with $k' \neq k$, so 
	the claim follows from our assumption that $k + k' \not\equiv 2 \bmod p-1$.  
	 
	Note that the condition
	\begin{equation} \label{cpcond}
		(\eta_{p-k},\eta_{k+k'-1})_{p,F,S} \neq 0
	\end{equation} 
	for all irregular $k' \neq k$ is equivalent to the statement that
	$$
		A_F^{(1-k')} \otimes \mu_p \subseteq 
		(\eta_{p-k},\mc{C}_F)_{p,F,S},
	$$
	as the odd eigenspaces of $A_F$ are cyclic and $\mc{C}_F^{(j)}$ is generated by
	$\eta_{p-j}$ for each even $j$.  Let $E$ denote the extension of $F$
	by a $p$th root of $\eta_{p-k}$.  Since $E/F$ is unramified, we have 
	$$
		A_{E/F} \cong \bigoplus_{\substack{k' \neq k\\
		(p,k') \text{ irregular}}} A_F^{(1-k')}.
	$$
	We therefore have 
	$$
		(\eta_{p-k},\mc{C}_F)_{p,F,S} = A_{E/F} \otimes \mu_p
	$$
	exactly when \eqref{cpcond} holds for any irregular $k'$ for $p$ 
	not equal to $k$. 
	
	If \eqref{cpcond} holds for all irregular $k' \neq k$,
	then Proposition \ref{surjtors} yields that $Y_L \cong X_{L/K}$
	via restriction, i.e., that $(I_H Y_L)_H = 0$ and, therefore, 
	$(I_H X_L)_H = 0$ by \eqref{samequot}.  On the other hand, if
	\eqref{cpcond} fails for some $k'$, then Lemma \ref{augtriv}
	implies that $I_H Y_L \neq 0$ and therefore that
	the restriction map $Y_L \to X_K$ is not injective, so in
	particular the restriction map $X_L \to X_K$ is not.  Therefore, in 
	this case, $\mf{G}_K$ is not abelian.
	
	We now suppose that \eqref{cpcond} holds for all $k$ and $k'$ irregular
	for $p$ with $k' > k$.  We have shown that $(I_H X_L)_H = 0$ for the field
	$L$ and Galois group $H$ determined by $k$ above, for each $k$.
	Since $H \wedge_{\zp} H = 0$, we have 
	by Lemma \ref{keyseq} that $X_L$ injects into $X_K$ with cokernel $H$.  
	Now let $M$ be the maximal unramified abelian 
	pro-$p$ extension of $K$, and set $G = X_K = \Gal(M/K)$
	and $N = X_L \cong \ker(G \to H)$.  
	From Lemma \ref{keyseq}, we obtain a commutative diagram with exact rows
	$$
		\SelectTips{cm}{} \xymatrix{
		N \wedge_{\zp} N \ar[d] \ar[r] & (X_M)_N \ar[r] \ar[d] & 
		X_L \ar[r]^{\sim} \ar[d] & N \ar[r] \ar[d] & 0 \\ 
		G \wedge_{\zp} G \ar[r] & (X_M)_G \ar[r] & X_K \ar[r]^{\sim} & 
		G \ar[r] & 0.
		}
	$$
	In the diagram, the leftmost horizontal arrows are surjective, the
	leftmost vertical arrow is the natural injection, and the second vertical
	arrow is the natural surjection.  It follows that $N \wedge_{\zp} N 
	\to (X_M)_G$ is surjective.  Now $G \wedge_{\zp} G$ (resp., 
	$N \wedge_{\zp} N$) is nontrivial only
	in its $(\omega^{2-j-j'})$-eigenspaces, where $j < j'$ are irregular for $p$
	(resp., are irregular for $p$ and distinct from $k$).
	Since the $j+j' \pmod{p-1}$ with $j < j'$ are all distinct by assumption, 
	we have that $(X_M)_G$ is trivial in its $\omega^{2-k-k'}$-eigenspaces 
	for all $k'$
	irregular for $p$ and distinct from $k$.  Since this holds for all $k$,
	we have $(X_M)_G = 0$.	This forces 
	$[\mf{G}_K,\mf{G}_K]^{\mr{ab}} \cong X_M = 0$, so $\mf{G}_K$ is abelian.
\end{proof}	
	
\begin{corollary}
	For $p < 1000$, the group $\mf{G}_K$ is abelian.  For $p = 1217$, $7069$,
	and $9829$, it is nonabelian.  For all other $p < 25,\!000$, it is abelian
	if $(\,\cdot\,,\cdot\,)_{p,F,S}$ is surjective.
\end{corollary}	

\begin{proof}
	Vandiver's conjecture and the assumption that the nontrivial eigenspaces of 
	$X_K$ are procyclic are known for $p < $\ 12,000,000 \cite{bcems}.
	Using Magma, we checked that the congruences assumed not to hold among 
	irregular $k$ for $p$
	are never satisfied for $p < $\ 25,000.
	We then checked the pairing values that appear in the tables mentioned 
	in \cite[Theorem 5.1]{mcs} (now completed for 
	$p < $\ 25,000).  That is,
	for each irregular $k$ and $k'$ for $p$ with $k < k'$, 
	we checked the table value $b_{p,k,k'}
	\in \Z/p\Z$ 
	that corresponds to $a_{p,k,k'} = (\eta_{p-k},\eta_{k+k'-1})_{p,F,S}$ 
	to see if it
	is zero.  Since 
	$$
		a_{p,k,k'} = b_{p,k,k'} \cdot \mf{c}_{p,k'}
	$$ 
	for some $\mf{c}_{p,k'}$ in 
	$A_F^{(1-k')} \otimes \mu_p$, the triviality of 
	$b_{p,k,k'}$ implies the triviality
	of $a_{p,k,k'}$.  On the other hand, the converse holds whenever
	$(\,\cdot\,,\,\cdot\,)_{p,F,S}$ is surjective.
	We know this surjectivity to be the case for $p < 1000$ by 
	Theorem \ref{pairp}. 
	
	The values $b_{p,k,k'}$ with $p < $\ 25,000 are nonzero aside from three
	exceptions: 
	\begin{itemize}
		\item[$\bullet$] $p = 1217$, $k = 784$, and $k' = 866$,
		\item[$\bullet$] $p = 7069$, $k = 1478$, and $k' = 2570$,
		\item[$\bullet$] $p = 9829$, $k = 4562$, and $k' = 7548$.
	\end{itemize}
	(We remark that the index of irregularity of $p$ is $3$ for $p = 1217$ and
	$2$ for $p = 7069$ and $p = 9829$.)  The result now follows from
	Theorem \ref{GKresult}.
\end{proof}

We remark that McCallum and the author have conjectured that $(\,\cdot\,,
\,\cdot)_{p,F,S}$ is surjective when Vandiver's conjecture holds at $p$
\cite[Conjecture 5.3]{mcs}, so we expect $\mf{G}_K$ to be abelian for
all $p < $\ 25,000 other than $1217$, $7069$, and $9829$.

\renewcommand{\baselinestretch}{1}

\vspace{2ex} \footnotesize \noindent
Romyar Sharifi \\
Department of Mathematics and Statistics\\
McMaster University\\
1280 Main Street West\\
Hamilton, Ontario L8S 4K1, Canada\\
e-mail address: {\tt sharifi@math.mcmaster.ca}\\
web page: {\tt http://www.math.mcmaster.ca/$\thicksim$sharifi}
\end{document}